\documentclass{amsart}

\usepackage{amsmath,amsthm,amsfonts,amssymb}

\usepackage{hyperref}
\hypersetup{colorlinks=true,citecolor=blue,linkcolor=blue,urlcolor=blue,
pdfstartview=FitH, pdfauthor=Gergely
Harcos,pdftitle=Twisted Hilbert modular L-functions and spectral theory}

\theoremstyle{plain}
\newtheorem{theorem}{Theorem}

\theoremstyle{definition}
\newtheorem*{acknowledgements}{Acknowledgements}

\renewcommand{\leq}{\leqslant}

\DeclareMathOperator{\GL}{GL} \DeclareMathOperator{\SL}{SL}
\DeclareMathOperator{\SO}{SO} 
 \DeclareMathOperator{\disc}{disc}
\DeclareMathOperator{\ord}{ord} \DeclareMathOperator{\lcm}{lcm}
 \DeclareMathOperator*{\aut}{aut}
\DeclareMathOperator*{\gen}{gen} \DeclareMathOperator*{\spn}{spn}

\newcommand{\ov}[1]{\overline{#1}}

\newcommand{\eps}{\varepsilon}

\newcommand{\RR}{\mathbb{R}}
\newcommand{\CC}{\mathbb{C}}
\newcommand{\QQ}{\mathbb{Q}}
\newcommand{\ZZ}{\mathbb{Z}}
\newcommand{\aA}{\mathbb{A}}

\newcommand{\NN}{\mathbb{N}}
\newcommand{\ma}{\mathfrak{a}}

\newcommand{\mc}{\mathfrak{c}}
\newcommand{\mo}{\mathfrak{o}}
\newcommand{\mpr}{\mathfrak{p}}
\newcommand{\mt}{\mathfrak{t}}
\newcommand{\mq}{\mathfrak{q}}
\newcommand{\my}{\mathfrak{y}}

\newcommand{\KK}{\mathcal{K}}
\newcommand{\N}{\mathcal{N}}
\newcommand{\Kplus}{K_{\infty,+}^{\times}}
\newcommand{\Kplusdiag}{K_{\infty,+}^{\text{diag}}}

\newcommand{\Afin}{\mathbb{A}_{\text{fin}}}
\newcommand{\mm}{\mathfrak{m}}
\newcommand{\cS}{\mathcal{S}}
\newcommand{\nyil}{\ \rightsquigarrow\ }

\begin{document}

\author{Gergely Harcos}
\address{Alfr\'ed R\'enyi Institute of Mathematics, Hungarian Academy of Sciences, POB 127, Budapest H-1364, Hungary}\email{gharcos@renyi.hu}
\address{Central European University, Nador u. 9, Budapest H-1051, Hungary}\email{harcosg@ceu.hu}

\title{Twisted Hilbert modular $L$-functions and spectral theory}

\thanks{The author was supported by OTKA grants K~101855 and K~104183 and ERC Advanced Grant~228005.}

\begin{abstract}
These are notes for four lectures given at the 2010 CIMPA Research School
``Automorphic Forms and $L$-functions'' in Weihai, China. The lectures focused on a Burgess-like subconvex bound for twisted Hilbert modular $L$-functions published jointly with Valentin Blomer in the same year. They discussed the proof in some detail, especially how spectral theory can be used to estimate the relevant shifted convolution sums efficiently. They also discussed briefly an application for the number of representations by a totally positive ternary quadratic form over a totally real number field. The notes below follow the leisurely style of the lectures, hence they do not constitute a comprehensive survey of the subject.
\end{abstract}

\maketitle

\section{Lecture One: Some Quadratic Forms}

In the lectures we shall discuss a state-of-the-art bound for twisted Hilbert modular $L$-functions, namely an analogue of Burgess' bound~\cite{Bu} for Dirichlet $L$-functions. For motivation and context, we start with an application to quadratic forms.

For a positive integral ternary quadratic form $Q$ let us denote by $r^*(n,Q)$ the number of \emph{primitive} integral representations of $n$ by $Q$, that is
\[r^*(n,Q):=\#\left\{(x,y,z)\in\ZZ^3:\text{ $Q(x,y,z)=n$ and $\gcd(x,y,z)=1$}\right\}.\]
Geometrically, this is the number of \emph{visible} lattice points on an ellipsoid of size $\sqrt{n}$. We would like to understand this quantity, in particular determine which positive integers $n$ are primitively represented by $Q$ over $\ZZ$. Note that when $n$ is square-free, all
integral representations are primitive.

\subsection{Sums of three squares}
Let us look at the simplest example,
\[Q(x,y,z):=x^2+y^2+z^2.\]
Clearly, for $n\equiv 0,4,7\pmod{8}$ there are no primitive representations. However, for $n\equiv 1,2,3,5,6\pmod{8}$ there are, and we have the following elegant formula for their number:
\begin{equation}\label{eq1}
r^*(n,Q)=\frac{24}{\pi}\,\sqrt{n}\,L\left(1,\left(\frac{D}{\cdot}\right)\right),
\end{equation}
where
\[D=\begin{cases}
-4n,&n\equiv 1,2,5,6\pmod{8},\\
-n,&n\equiv 3\pmod{8}.
\end{cases}\]
This formula arises naturally from the work of Gauss~(1801) and Dirichlet~(1839) on the class number.
For a negative discriminant $D$ let $h(D)$ denote the number of equivalence classes of
positive \emph{primitive} binary quadratic forms of discriminant $D$, and let $w$ denote the number of automorphs\footnote{In these notes, equivalence classes and automorphs are meant in the narrow/strict/proper sense, i.e. they are defined in terms of the action of $\SL_n$ on $n$-ary quadratic forms.} of such a form:
\[w=\begin{cases}6,&D=-3,\\4,&D=-4,\\2,&D<-4.\end{cases}\]
In this notation a fundamental result of Gauss~\cite{Ga} states that
\begin{equation}\label{eq2}
r^*(n,Q)=\begin{cases}
\frac{24}{w}h(-4n),&n\equiv 1,2,5,6\pmod{8},\\
\frac{48}{w}h(-n),&n\equiv 3\pmod{8},
\end{cases}
\end{equation}
while the class number formula of Dirichlet~\cite{D} reads
\begin{equation}\label{eq3}
h(D)=\frac{w}{2\pi}\,|D|^{1/2}\,L\left(1,\left(\frac{D}{\cdot}\right)\right).
\end{equation}
For a modern account of these results we refer the reader to \cite[\S 4]{G} and \cite[\S 8]{Z}, including the exercises. Combining \eqref{eq2} and \eqref{eq3}, we obtain \eqref{eq1}.

It is not hard to see that the special $L$-value involved in \eqref{eq1}
is not too large, while a deeper result of Siegel~\cite{S2} shows that it is not too small either:
\begin{equation}\label{eq4}
L\left(1,\left(\frac{D}{\cdot}\right)\right)=|D|^{o(1)}.
\end{equation}
We infer that $r^*(n,Q)$ is either zero or about $\sqrt{n}$:
\begin{equation}\label{eq5}
r^*(n,Q)=n^{\frac{1}{2}+o(1)},\qquad n\equiv 1,2,3,5,6\pmod{8}.
\end{equation}
That is, for the case of $Q(x,y,z)=x^2+y^2+z^2$, we know precisely when there is at least one visible lattice point on our ellipsoid (the sphere of radius $\sqrt{n}$ centered at the origin), and we also know that one visible lattice point implies many others.
We can look at \eqref{eq5} as a quantitative local-to-global principle:
\[\boxed{\text{$n$ is primitively represented by $Q$ over any $\ZZ_p$}\quad\Longrightarrow\quad\text{$r^*(n,Q)$ is large}}\]
Here we could restrict the left hand side to $p=2$, because for $p>2$ the condition is automatically met.
This formulation, of course, does not reveal that $r^*(n,Q)$ has a nice analytic description as in \eqref{eq1}.
The two viewpoints are unified by the celebrated mass formula\footnote{The original formula concerns all integral representations, but it is easy to reformulate it in terms of primitive representations for the cases considered here. For the general case see \cite[\S 6.8]{K}.} of Siegel~\cite{S1}. Namely, let us assume that $n\equiv 1,2,3,5,6\pmod{8}$ so that $x^2+y^2+z^2=n$ has a solution in $\RR$ and a primitive solution in each $\ZZ_p$. Then we can calculate the right hand side of \eqref{eq1} as a product of densities of the primitive local solutions over $\RR$ and the various completions $\ZZ_p$:
\[\beta_\infty=2\pi\sqrt{n},\qquad\beta_2=\frac{\frac{3}{2}}{1-(\frac{D}{2})\frac{1}{2}},
\qquad\beta_p=\frac{1-\frac{1}{p^2}}{1-(\frac{D}{p})\frac{1}{p}}\quad (p\neq 2),\]
whose product is indeed
\[\beta_\infty\prod_p\beta_p=2\pi\sqrt{n}\,\frac{\frac{3}{2}}{1-\frac{1}{2^2}}\,\frac{1}{\zeta(2)}\,L\left(1,\left(\frac{D}{\cdot}\right)\right)=
\frac{24}{\pi}\,\sqrt{n}\,L\left(1,\left(\frac{D}{\cdot}\right)\right).\]

\subsection{Ramanujan's ternary quadratic form}
What about more general ternary quadratic forms? As a second example let us consider\footnote{This section was inspired by a nice paper of Ono and Soundararajan~\cite{OS}.}
\[Q(x,y,z):=x^2+y^2+10z^2.\]
Is there a similar elegant formula as \eqref{eq1} for the number of primitive representations?
Well, almost. What happens is that there is another ternary quadratic form which is equivalent to
$Q$ over $\RR$ and all the completions $\ZZ_p$, but not over $\ZZ$:
\[Q'(x,y,z):=2x^2+2y^2+3z^2-2xz.\]
So $Q$ and $Q'$ produce the same primitive local densities of representations
over the reals and the $p$-adic integers, yet they are
really different over the integers. A quadratic form that is locally equivalent to $Q$ is globally equivalent to $Q$ or $Q'$,
hence perhaps it is not surprising that the
product of primitive local densities is related to a combination of $r^*(n,Q)$ and $r^*(n,Q')$. For this case Siegel's
mass formula~\cite{S1} precisely tells us that
\begin{equation}\label{eq6}
\frac{1}{3}r^*(n,Q)+\frac{2}{3}r^*(n,Q')=\frac{4\sqrt{10}}{3\pi}\,\sqrt{n}\,L\left(1,\left(\frac{-10n}{\cdot}\right)\right),
\end{equation}
at least when $\gcd(n,10)=1$ as we assume for simplicity. Here the weight of $r^*(n,Q)$ is one-half of the weight of
$r^*(n,Q')$, because $Q$ has twice as many automorphs as $Q'$ (eight vs. four). The right hand side manifests as the product of the
common primitive local densities of $Q$ and $Q'$,
\[\beta_\infty=\frac{2\pi\sqrt{n}}{\sqrt{10}},\qquad\beta_2=1,\qquad\beta_5=\frac{4}{5},
\qquad\beta_p=\frac{1-\frac{1}{p^2}}{1-(\frac{-10n}{p})\frac{1}{p}}\quad (p\neq 2,5).\]
Indeed, the product of these quantities equals
\[\frac{2\pi\sqrt{n}}{\sqrt{10}}\ \frac{1}{1-\frac{1}{2^2}}\,\frac{\frac{4}{5}}{1-\frac{1}{5^2}}\,\frac{1}{\zeta(2)}\,L\left(1,\left(\frac{-10n}{\cdot}\right)\right)=
\frac{4\sqrt{10}}{3\pi}\,\sqrt{n}\,L\left(1,\left(\frac{-10n}{\cdot}\right)\right).\]

We are still assuming that $\gcd(n,10)=1$. Comparing \eqref{eq6} with \eqref{eq4}, we see that $n$ has many primitive representations by $Q$ or $Q'$. However, we would like to know if $n$ has any or many primitive representations by $Q$ alone! Using automorphic forms and $L$-functions one can show that $r^*(n,Q)\approx r^*(n,Q')$ with great precision, so that both $r^*(n,Q)$ and $r^*(n,Q')$ are close to the right hand side of \eqref{eq6}. Specifically, by the work of Schulze-Pillot~\cite{SP1} and Waldspurger~\cite{W} based on the Shimura correspondence~\cite{Sh}, there exists a primitive holomorphic cusp form $f$ of weight $2$, level $20$, and trivial nebentypus such that
\begin{equation}\label{eq7}
\bigl(r^*(n,Q)-r^*(n,Q')\bigr)^2\ll_\eps n^{\frac{1}{2}+\eps}\,L\left(\frac{1}{2},f\otimes\biggl(\frac{-10 n}{\cdot}\biggr)\right).
\end{equation}

If we assume the Grand Riemann Hypothesis (GRH) for automorphic $L$-functions, then the central $L$-value on the right hand side of \eqref{eq7} is $\ll_\eps n^\eps$, whence from \eqref{eq6} we infer that
\begin{equation}\label{eq8}
r^*(n,Q)=\frac{4\sqrt{10}}{3\pi}\,\sqrt{n}\,L\left(1,\left(\frac{-10n}{\cdot}\right)\right)+O_\eps\left(n^{\frac{1}{4}+\eps}\right).
\end{equation}
The main term here is $\gg_\eps n^{\frac{1}{2}-\eps}$ by \eqref{eq4}, and all the implied constants are effective under GRH. Roughly, this means that
\[\boxed{\text{GRH and $\gcd(n,10)=1$}\quad\Longrightarrow\quad\text{$r^*(n,Q)$ is large}}\]
This is another instance (albeit conditionally) of a quantitative local-to-global principle, because the integers coprime with $10$ are primitively represented by $Q$ over any $\ZZ_p$. In particular, under GRH, there are only finitely many positive integers such that $\gcd(n,10)=1$ and $r^*(n,Q)=0$, and GRH even provides a \emph{theoretical} algorithm to determine them. It is a challenging task to turn GRH into a \emph{practical} algorithm for this problem, solved for $n$ square-free by Ono and Soundararajan~\cite{OS}.

We can also estimate the right hand side of \eqref{eq7} without GRH. The functional equation for the twisted $L$-function coupled with the Phragm\'en-Lindel\"of convexity principle shows that the $L$-value in \eqref{eq7} is $\ll_\eps n^{\frac{1}{2}+\eps}$. This \emph{convexity bound}, however, yields the error term $O_\eps\bigl(n^{\frac{1}{2}+\eps}\bigr)$ in \eqref{eq8}, which is too weak to imply $r^*(n,Q)>0$. If, instead, we employ the \emph{subconvexity bound}
\[L\left(\frac{1}{2},f\otimes\biggl(\frac{-10 n}{\cdot}\biggr)\right)\ll_\eps n^{\frac{1}{2}-\delta+\eps}\]
for some $\delta>0$, then we arrive at a useful variant of \eqref{eq8} for $\gcd(n,10)=1$:
\[r^*(n,Q)=\frac{4\sqrt{10}}{3\pi}\,\sqrt{n}\,L\left(1,\left(\frac{-10n}{\cdot}\right)\right)+O_\eps\left(n^{\frac{1-\delta}{2}+\eps}\right).\]
Such a conclusion was first established unconditionally with $\delta=\frac{1}{14}$ by Duke~\cite{Du} using a method of Iwaniec~\cite{Iw}
(see also \cite{DSP}), and then with $\delta=\frac{1}{8}$ by Bykovski\u\i~\cite{By} using a method of Duke, Friedlander, Iwaniec~\cite{DFI} (see also \cite{BH1} and \cite{HH}).
Furthermore, it seems likely that the technical assumptions in a deep result of Conrey and Iwaniec~\cite{CI} can be relaxed so as to yield the
conclusion with $\delta=\frac{1}{6}$. To summarize,
\[\boxed{\text{Subconvexity and $\gcd(n,10)=1$}\quad\Longrightarrow\quad\text{$r^*(n,Q)$ is large}}\]

\medskip

\section{Lecture Two: More Quadratic Forms}

The examples presented in the first lecture can be generalized to a large extent. Let $K$ be a totally real number field with discriminant $D$ and
ring of integers $\mo$. Let $(a_{ij})$ be a $3\times 3$ matrix with $a_{ij}=a_{ji}\in\mo$ and $a_{ii}\in 2\mo$ such that the corresponding integral quadratic form
\[Q(x_1,x_2,x_3):=\frac{1}{2}\sum_{i,j}a_{ij}x_ix_j\]
is totally positive. Then the determinant $d:=\det(a_{ij})$ is totally positive and lies in $2\mo$ by \cite[Lemma~2.1]{H}.
We are interested in the number of \emph{primitive} integral representations of a totally positive integer $n\in\mo$ by $Q$, that is
\[r^*(n,Q):=\#\left\{(x,y,z)\in\mo^3:\text{ $Q(x,y,z)=n$ and $\gcd(x,y,z)=\mo$}\right\}.\]

An obvious necessary condition for $r^*(n,Q)>0$ is that $n$ is primitively represented by $Q$ over any
non-archimedean completion $\mo_\mpr$ of $\mo$. This condition is invariant under replacing $Q$ by any form in its genus, that is, by any totally positive ternary quadratic form over $\mo$ which is locally equivalent to $Q$ over any $\mo_\mpr$. In fact the primitive local densities of the representations only depend on the genus, and by Siegel's mass formula~\cite{S1} their product equals a weighted average
over the finitely many equivalence classes contained in the genus:
\[r^*(n,\gen Q):=\left(\sum_{[Q']\in\gen Q}\frac{r^*(n,Q')}{\aut(Q')}\right)\left(\sum_{[Q']\in\gen Q}\frac{1}{\aut(Q')}\right)^{-1} = \beta_\infty\prod_\mpr\beta_\mpr.\]
Here $r^*(n,Q')$ and $\aut(Q')$, the number of automorphs of $Q'$, only depend on the class $[Q']$ of $Q'$, hence the sums are well-defined.
For $\gcd(n,d)=\mo$ we can simplify the right hand side using \cite[\S 7 of Part III]{S1} and \cite[Lemma~3.2]{H} to obtain
\begin{equation}\label{eq9}
r^*(n,\gen Q)=c(n,Q)\,(\N n)^{\frac{1}{2}}\,L\left(1,\left(\frac{-2dn}{\cdot}\right)\right),
\end{equation}
where $\N$ stands for the norm and $c(n,Q)$ equals, up to a positive constant depending on $K$, the density of primitive solutions of the congruence $Q(x,y,z)\equiv n\pmod{4d}$. In particular, $c(n,Q)>0$ is equivalent to primitive local representability modulo~$4d$, in which case $r^*(n,\gen Q)=(\N n)^{\frac{1}{2}+o(1)}$ by a straightforward extension of \eqref{eq4}.

Surprisingly, $r^*(n,\gen Q)$ being large does not ensure that $n$ is primitively represented by $Q$, even when $\N n$ is large and coprime with $d$. Here are two examples, borrowed from \cite{EHH} and \cite{SP2}. When $K=\QQ$, the forms $x^2+3y^2+36z^2$ and $3x^2+4y^2+9z^2$ are in the same genus, yet any square number coprime with $6$ is primitively represented by exactly one of them. When $K=\QQ(\sqrt{35})$, an integer of the form $7p^2$ with a rational prime $p\equiv 1\pmod{7}$ is not a sum of three squares in $\mo$,
although a sum of three coprime squares in $\mo_\mpr$ for any prime ideal $\mpr$. The proper discussion of this phenomenon would lead us too far as it relies on the theory of spinor genera and theta series, see the recent surveys \cite{H2,SP2} and the references therein. Let us just say that in the modern theory one considers lattices and their representations in a quadratic space over $K$, and the restriction to free lattices yields the classical theory of integral quadratic forms and their representations \cite[\S 82]{M}. A class is an orbit of lattices under the group of rotations of the space. A genus is an orbit under the group of adelic rotations, while a spinor genus is an orbit under a certain normal subgroup of adelic rotations \cite[\S 102]{M}. Lattices in the same genus are isomorphic as $\mo$-modules\footnote{This was kindly explained to me by Rainer Schulze-Pillot, here is a variant of his argument. Let $L$ be a lattice in a fixed quadratic space $V$ over $K$. By \cite[81:3]{M}, there is a basis $(x_i)$ of $V$ and fractional ideals $\ma_i$ in $K$ such that $L=\ma_1x_1+\dots+\ma_nx_n$. Writing $\ma=\ma_1\dots \ma_n$, the volume of $L$ equals $\ma^2\disc(x_1,\dots,x_r)$, which is $\ma^2\disc(V)$ modulo $(F^\times)^2$. Therefore the volume of $L$ (hence also the genus of $L$) determines $\ma$ modulo $F^\times$, which is the Steinitz class of $L$ as an $\mo$-module.}, hence the set of free lattices is a disjoint union of genera (resp. spinor genera) in the modern sense.

We shall conveniently avoid the above type of exceptions by restricting to totally positive \emph{non-square} integers $n\in\mo$ coprime with $dD$. Let $r^*(n,\spn Q)$ be defined similarly as $r^*(n,\gen Q)$, but over the spinor genus of $Q$. By the work of Kneser~\cite{Kn} and Hsia~\cite{Hs}, the two averages agree now\footnote{Let $L$ be a free ternary lattice corresponding to the class of $Q$. It suffices to show that $r(n,\spn\ma L)=r(n,\gen\ma L)$ for any ideal $\ma$ in $\mo$. Take a prime ideal $\mpr\nmid dD$ such that $\ord_\mpr(n)$ is odd. The quadratic extension $K_\mpr(\sqrt{-2dn})/K_\mpr$ is ramified, while $L_\mpr$ is unimodular, hence local class field theory combined with \cite[92:5]{M} shows that \cite[(4)]{SP4} fails at $\mpr$ in the present setting.}, so that the bounds usually proved for $r^*(n,Q)-r^*(n,\spn Q)$ apply for $r^*(n,Q)-r^*(n,\gen Q)$. By the work of Baruch--Mao~\cite{BM}, Blasius~\cite{Bl}, Schulze-Pillot~\cite{SP3}, Waldspurger~\cite{W2}, there exist finitely many primitive holomorphic Hilbert cusp forms $f_1,\dots, f_r$ over $K$ depending on $Q$, each of weight $(2,\dots,2)$ and trivial nebentypus, such that
\begin{equation}\label{eq10}
\bigl(r^*(n,Q)-r^*(n,\gen Q)\bigr)^2\ll_{Q,K,\eps} (\N n)^{\frac{1}{2}+\eps}\,
\max_{1\leq i\leq r}L\left(\frac{1}{2},f_i\otimes\biggl(\frac{-2dn}{\cdot}\biggr)\right).
\end{equation}

As before, any subconvex bound for the central $L$-values on the right hand side yields an asymptotic formula for $r^*(n,Q)$ with a power saving error term. More specifically and generally, let us assume the bound
\begin{equation}\label{eq11}
L\left(\frac{1}{2},\pi\otimes \chi\right) \ll_{\pi,\chi_\infty,K,\eps} (\N\mq)^{\frac{1}{2}-\delta+\eps},
\end{equation}
where $\delta$ is a positive constant, $\pi$ is any irreducible unitary cuspidal representation of $\GL_2$ over $K$, and $\chi$ is any Hecke character of $K$ of conductor $\mq$. Then from \eqref{eq9} and \eqref{eq10} we infer, for totally positive non-square $n\in\mo$ coprime with $dD$,
\begin{equation}\label{eq12}
r^*(n,Q)=c(n,Q)\,(\N n)^{\frac{1}{2}}\,L\left(1,\left(\frac{-2dn}{\cdot}\right)\right)+O_{Q,K,\eps}\left((\N n)^{\frac{1-\delta}{2}+\eps}\right).
\end{equation}
Recall that $c(n,Q)>0$ is equivalent to primitive local representability of $n$ by $Q$ modulo~$4d$, in which case the main term dominates the error term for $\N n$ sufficiently large:
\[\boxed{\text{$\gcd(n,dD)=\mo$ and $n\neq\square$ and $c(n,Q)>0$}\quad\Longrightarrow\quad\text{$r^*(n,Q)$ is large}}\]

\medskip

The conclusion \eqref{eq11} was first established by Cogdell, Piatetski-Shapiro, Sarnak~\cite{CPSS} with the value $\delta=\frac{1-2\theta}{14+4\theta}$, at least for $\pi$ induced by a holomorphic Hilbert cusp form, which suffices for the application \eqref{eq12}. For general $\pi$ (and arbitrary $K$), the breakthrough is due to Venkatesh~\cite{Ve} who achieved $\delta=\frac{(1-2\theta)^2}{14-2\theta}$. In these results, $0\leq\theta\leq \frac{1}{2}$ is an approximation towards the Ramanujan--Petersson conjecture, the current record being $\theta=\frac{7}{64}$ due to Blomer and Brumley~\cite{BB}.
Interestingly, under the Ramanujan--Petersson conjecture $\theta=0$ both expressions become $\delta=\frac{1}{14}$, matching the already mentioned results of Duke~\cite{Du} and Iwaniec~\cite{Iw}. For totally real $K$ and general $\pi$, Blomer and Harcos~\cite{BH2} established the Burgess-like subconvexity saving $\delta=\frac{1-2\theta}{8}$, and the proof of this result is outlined in the rest of these notes. Recently Maga~\cite{Ma,Ma2} extended the method of \cite{BH2} to arbitrary number fields. Wu~\cite{Wu} obtained the same $\delta=\frac{1-2\theta}{8}$ in full generality, even uniformly in $\chi_\infty$, by a different method based on the deep work of Michel and Venkatesh~\cite{MV}.

\begin{theorem}[\cite{BH2}]\label{thm1} Let $K$ be a totally real number field. Let $\pi$ be an irreducible unitary cuspidal representation of $\GL_2$ over $K$, and $\chi$ a Hecke character of $K$ of conductor $\mq$. Then for any $\eps>0$ one has
\[L\left(\frac{1}{2},\pi\otimes \chi\right) \ll_{\pi,\chi_\infty,K,\eps} (\N\mq)^{\frac{3+2\theta}{8}+\eps}.\]
\end{theorem}

\begin{theorem}[\cite{BM,Bl,BH2,SP3,S1,W2}]\label{thm2} Let $K$ be a totally real number field with discriminant $D$. Let $Q$ be a totally positive integral ternary quadratic form over $K$ with determinant $d$. Let $n$ be a totally positive non-square integer in $K$ coprime to $dD$. Then the number of primitive representations $Q(x,y,z)=n$ equals
\[r^*(n,Q)=c(n,Q)\,(\N n)^{\frac{1}{2}}\,L\left(1,\left(\frac{-2dn}{\cdot}\right)\right)+O_{Q,K,\eps}\left((\N n)^{\frac{7+2\theta}{16}+\eps}\right),\]
where $c(n,Q)$ is a constant times the density of primitive solutions of the congruence $Q(x,y,z)\equiv n\pmod{4d}$.
\end{theorem}

\medskip

\section{Lecture Three: Preliminaries from Number Theory}

We collect some preliminaries for the proof of Theorem~\ref{thm1}, to be outlined in the fourth lecture. As explained in the first two lectures, Theorem~\ref{thm2} follows from Theorem~\ref{thm1}. For more detail concerning the preliminaries the reader should consult \cite{BH2,GJ,N,We}.

\subsection{Adeles and ideles}
The adele ring of $K$ is a restricted direct product of the completions $K_v$ at the various places of $K$:
\[\aA:=K_\infty \times \Afin=\prod_{v\mid\infty}K_v\times{\prod_\mpr}' K_\mpr,\]
where $\mpr$ runs through the prime ideals of $\mo$. The topology of the additive group $(\aA,+)$ is determined by the fundamental neighborhoods of zero in $\Afin$:
\[U^\mm:=\prod_\mpr U_\mpr^{(\ord_\mpr\mm)},\]
where $\mm\subseteq\mo$ is any nonzero ideal, $\ord_\mpr\mm$ denotes the exponent of $\mpr$ in $\mm$, and
\[U_\mpr^{(n)}:=\mpr^n\mo_\mpr,\qquad n\in\NN.\]
Then $(\aA,+)$ is a locally compact Hausdorff topological group such that
\[K\overset{\text{diag}}\hookrightarrow\aA\text{ is discrete},\qquad K\backslash\aA\text{ is compact}.\]
Similarly, the group of ideles of $K$ is a restricted direct product
\[\aA^\times:=K_\infty^\times \times \Afin^\times=\prod_{v\mid\infty}K_v^\times\times{\prod_\mpr}' K_\mpr^\times,\]
whose topology is determined by the fundamental neighborhoods of one in $\Afin^\times$:
\[V^\mm:=\prod_\mpr V_\mpr^{(\ord_\mpr\mm)},\]
where $\mm\subseteq\mo$ is any nonzero ideal, $\ord_\mpr\mm$ denotes the exponent of $\mpr$ in $\mm$, and
\[V_\mpr^{(n)}:=\begin{cases}\mo_\mpr^\times,&n=0;\\1+\mpr^n\mo_\mpr,&n>0.\end{cases}\]
Then $(\aA^\times,\cdot)$ is a locally compact Hausdorff topological group such that
\[\aA^\times=\Kplusdiag\aA^1\cong \RR_{>0}\times\aA^1,\qquad K^\times\overset{\text{diag}}\hookrightarrow\aA^1\text{ is discrete},\qquad
K^\times\backslash\aA^1\text{ is compact}.\]
Here $\aA^1$ is the subgroup of ideles of module $1$, and $\Kplusdiag$ denotes $\RR_{>0}$ diagonally embedded into $K_\infty^\times$.

\subsection{Hecke characters and Gr\"ossencharacters}
A \emph{Hecke character} is a continuous homomorphism $\chi:\aA^\times\to S^1$ which is trivial on $K^\times$. The kernel of such a character $\chi$ always contains a subgroup of the form $\{(1,\dots,1)\}\times V^\mm$. Let $I^\mm$ denote the group of fractional ideals of $K$ coprime with $\mm$, this is a free abelian group generated by the prime ideals $\mpr\nmid\mm$. Choose a prime element $p_\mpr\in\mpr-\mpr^2$ for each $\mpr\nmid\mm$, and define the character $\tilde\chi:I^\mm\to S^1$ via
\[\tilde\chi(\mpr):=\chi(\dots,1,1,p_\mpr,1,1,\dots),\qquad\mpr\nmid\mm.\]
Observe that $\tilde\chi$ is independent of the choice of the $p_\mpr$'s. Moreover, for any \emph{principal} ideal $(a)\in I^\mm$ have
\begin{align*}
\tilde\chi((a))
&=\chi\bigl(\underbrace{1,\dots,1}_{\text{arch}},\underbrace{\dots,p_\mpr^{\ord_\mpr(a)},\dots}_{\text{non-arch}}\bigr)\\
&=\chi\bigl(\underbrace{1,\dots,1}_{v\mid\infty},\underbrace{1,\dots,1}_{\mpr\mid\mm},\underbrace{a,a,\dots}_{\mpr\nmid\mm}\bigr)\\
&=\chi\bigl(\underbrace{a^{-1},\dots,a^{-1}}_{v\mid\infty},\underbrace{a^{-1},\dots,a^{-1}}_{\mpr\mid\mm},\underbrace{1,1,\dots}_{\mpr\nmid\mm}\bigr).
\end{align*}
If we interpret $a\in K^\times$ as an element of $K_\infty^\times$ by embedding $a$ diagonally, and as an element of $(\mo/\mm)^\times$ by reducing $a$ modulo~$\mm$, then we infer that
\[\tilde\chi((a))=\tilde\chi_\infty(a)\tilde\chi_\text{fin}(a),\qquad (a)\in I^\mm,\]
where $\tilde\chi_\infty:K_\infty^\times\to S^1$ and $\tilde\chi_\text{fin}:(\mo/\mm)^\times\to S^1$ are uniquely determined characters\footnote{In these notes, all characters are continuous.}. A character $\tilde\chi:I^\mm\to S^1$ with this property is called a \emph{Gr\"ossencharacter} modulo~$\mm$. So any Hecke character can be regarded as a Gr\"ossencharacter.

Conversely, any Gr\"ossencharacter $\tilde\chi:I^\mm\to S^1$ arises in this way from a Hecke character $\chi:\aA^\times\to S^1$ trivial on $\{(1,\dots,1)\}\times V^\mm$. To see this, define
\begin{align*}
V&:=K_\infty^\times\times\prod_{\mpr\mid\mm}V_\mpr^{(\ord_\mpr\mm)}\prod_{\mpr\nmid\mm}K_\mpr^\times;\\
\chi(a)&:=\tilde\chi_\infty^{-1}\left(a_\infty\right)\tilde\chi\bigl(\prod_{\mpr\nmid\mm}\mpr^{\ord_\mpr a_\mpr}\bigr),\qquad a\in V.
\end{align*}
Observe that $V$ is a subgroup of $\aA^\times$, and $\chi:V\to S^1$ is a character trivial on $\{(1,\dots,1)\}\times V^\mm$. This character extends uniquely to a character $\chi:K^\times V\to S^1$ trivial on $K^\times$, because for $a\in K^\times\cap V$ we have $a\equiv 1\pmod{\mm}$, whence
\[\chi(a)=\tilde\chi_\infty^{-1}(a)\tilde\chi((a))=\tilde\chi_\text{fin}(a)=1,\qquad a\in K^\times\cap V.\]
On the other hand, $K^\times V$ equals $\aA^\times$ by weak approximation in $K$, so we are done.

If, for given $\chi$, we choose the largest $\mm$ (called the \emph{conductor}), then $\tilde\chi$ will be primitive, and vice versa. To summarize, we have a natural bijection
\[\boxed{\text{$\chi$ a Hecke character}\quad\longleftrightarrow\quad\text{$\tilde\chi$ a primitive Gr\"ossencharacter}}\]
In the rest of these notes, $\chi$ will always stand for a Hecke character, and $\tilde\chi$ for the corresponding primitive Gr\"ossencharacter.

\subsection{The automorphic spectrum of $\GL_2$}
Let us restrict to trivial central character for simplicity. Then the corresponding $L^2$ space of automorphic functions on $\GL_2$ over $K$ decomposes as a direct integral of irreducible automorphic representations,
\[L^2(\GL_2(K)\backslash \GL_2(\mathbb{A}))=
\left(\bigoplus_{\pi} V_{\pi}\right)\oplus
\left(\bigoplus_{\chi^2=1} V_\chi\right)\oplus
\left(\int\limits_{\{\chi,\chi^{-1}\}}V_{\chi,\chi^{-1}}\,d\{\chi,\chi^{-1}\}\right),\]
in the sense that each function in the $L^2$ space is a convergent integral of
functions from each subspace, and a Plancherel formula holds. In this decomposition, which is compatible with the
right action of $\GL_2(\mathbb{A})$,
\begin{itemize}
\item $V_\pi$ is an irreducible subspace of the cuspidal space defined by
\[\int_{K\backslash\mathbb{A}}\phi\left(\left(\begin{matrix}1 & x\\0 & 1\end{matrix}\right) g\right) dx = 0 \quad\text{for almost all } g \in \GL_2(\mathbb{A}).\]\vskip 6pt
\item $V_\chi$ is the one-dimensional subspace spanned by the function $g\mapsto\chi(\det g)$.\vskip 8pt
\item $V_{\chi,\chi^{-1}}$ consists of the Eisenstein series
\[E(\varphi,g)\doteq \sum_{\gamma \in P(K)\backslash \GL_2(K)} \varphi(\gamma g),\]
where $P$ stands for upper triangular matrices, and $\varphi : \GL_2(\mathbb{A}) \to \CC$ satisfies
\[\int_{\KK} |\varphi(k)|^2 dk < \infty,\qquad \KK:=\SO_2(K_\infty)\times\KK(\mo);\]
\[\varphi\left(\begin{pmatrix} a & x\\0 & b\end{pmatrix}g\right) = \chi\left(\frac{a}{b}\right) \left|\frac{a}{b}\right|^{1/2}\varphi(g), \qquad
\begin{pmatrix}a & x\\0 & b\end{pmatrix}\in P(\aA).\]
Here $\KK(\mo)$ is a maximal compact subgroup of $\GL_2(\Afin)$ to be defined in the next subsection.
More precisely, the sum converges only when the exponent $1/2$ is replaced by any $s\in\CC$ with $\Re s>1$. The symbol $\doteq$ stands for evaluation at $s=1/2$ of the function obtained by meromorphic continuation to $s\in\CC$, keeping the restriction of $\varphi$ to $\KK$
fixed all the way.
\end{itemize}

\subsection{Conductor, $L$-function, Kirillov model}
Let $\mathfrak{d}$ denote the different ideal of $K$. For each cuspidal space $V_\pi$ (and also for each Eisenstein space $V_{\chi,\chi^{-1}}$) there is a largest congruence subgroup
\[\KK(\mm):=\prod_\mpr\biggl\{\begin{pmatrix}a & b\\c & d\end{pmatrix}\in\GL_2(K_\mpr):\ a,d \in \mo_\mpr,\ b \in \mathfrak{d}_\mpr^{-1}, \ c \in \mm\mathfrak{d}_\mpr,\ ad-bc\in\mo_\mpr^\times\biggr\},\]
regarded as a subgroup of $\GL_2(\aA)$ in the obvious way, such that $V_\pi$ contains a nonzero vector fixed by the right action of $\KK(\mm)$. The corresponding nonzero ideal $\mm$ is called the \emph{conductor} of $\pi$, denoted $\mc_\pi$. It is known that $\mc_{\chi,\chi^{-1}}=\mc_\chi^2$. If $V_\pi(\mm)$ denotes the $\KK(\mm)$-fixed subspace of $V_\pi$, then $V_\pi(\mc_\pi)$ (the space of \emph{newforms}) is particularly nice. Namely,
for each nonzero ideal $\mm$, the Hecke operator $T(\mm)$ acts by some scalar $\lambda_{\pi}(\mathfrak{m})\in\CC$ on this space. These Hecke eigenvalues satisfy $\lambda_{\pi}(\mathfrak{m})\ll_{\eps}(\N\mm)^{\theta+\eps}$ with $\theta=\frac{7}{64}$ by~\cite{BB}, and they
determine the $L$-function of~$\pi$:
\begin{align*}
L(s,\pi)&=\sum_{\{0 \}\neq \mathfrak{m} \subseteq \mo} \frac{\lambda_{\pi}(\mathfrak{m})}{(\N\mathfrak{m})^{s}}\\
&=\prod_{\mpr\mid\mc_\pi}\frac{1}{1-\lambda_\pi(\mpr)(\N\mpr)^{-s}}\prod_{\mpr\nmid\mc_\pi}\frac{1}{1-\lambda_\pi(\mpr)(\N\mpr)^{-s}+(\N\mpr)^{-2s}}.
\end{align*}
Fixing a nontrivial character $\psi:K\backslash\aA\to S^1$, any newform $\phi\in V_\pi(\mc_\pi)$ has a Fourier expansion
\[\phi\left(\left(\begin{matrix} y & x\\ 0 &
1\end{matrix}\right)\right) = \sum_{r \in K^\times}
\frac{\lambda_{\pi}(r y_{\text{fin}})}{\sqrt{\N (r y_{\text{fin}})}}
W_{\phi}(r y_{\infty}) \psi(rx),\qquad y\in \mathbb{A}^{\times},\ \
x \in \mathbb{A},\]
where $W_\phi\in L^2(K_\infty^\times,d^\times y)$ is given by
\[W_{\phi}(y) := \int_{K\backslash\mathbb{A}} \phi\left(\begin{pmatrix} y & x\\ 0 &
1\end{pmatrix}\right) \psi(-x) \,dx,\qquad y\in K_\infty^\times.\]
So we restricted here the Whittaker model to upper triangular matrices (this is called the Kirillov model), and
we separated the archimedean and non-archimedan parts. An important feature is that any Whittaker function $W_\phi\in L^2(K_\infty^\times,d^\times y)$ occurs for some newform $\phi\in V_\pi(\mc_\pi)$, and $\|W_\phi\|$ is proportional to $\|\phi\|$ with a constant depending very mildly on~$\pi$.

What about ``oldforms of level $\mc$'', i.e.\ the elements of $V_\pi(\mc)$ for $\mc$ divisible by $\mc_\pi$? For any nonzero ideal $\mt\mid\mc\mc_\pi^{-1}$ we have an isometric embedding of complex vector spaces
\[R_\mt:V_\pi(\mc_\pi)\hookrightarrow V_\pi(\mc),\qquad
(R_\mt\phi)(g):= \phi\left(g\begin{pmatrix}t^{-1}&0\\0&1\end{pmatrix}\right),\]
where $t\in\mathbb{A}_\text{fin}^\times$ is any finite idele representing $\mt$. Then
it follows from the local result of Casselman~\cite{Ca} or the global result of Miyake~\cite{Mi} that
\[V_\pi(\mc)=\bigoplus_{\mt\mid\mc\mc_{\pi}^{-1}} R_\mt V_\pi(\mc_\pi),\qquad \mc_\pi\mid\mc.\]
A technical difficulty here is that the spaces $R_\mt V_\pi(\mc_\pi)$ are in general not orthogonal. Nevertheless, using a Gram--Schmidt
orthogonalization process, we can obtain an orthogonal decomposition
\[V_\pi(\mc)=\bigoplus_{\mt\mid\mc\mc_{\pi}^{-1}}R^{(\mt)}V_\pi(\mc_\pi),\qquad \mc_\pi\mid\mc,\]
and for every $\phi\in R_\mt V_\pi(\mc_\pi)$ a Fourier expansion
\[\phi\left(\left(\begin{matrix} y & x\\ 0 & 1\end{matrix}\right)\right) =
\sum_{r \in K^{\times}} \frac{\lambda^{(\mt)}_{\pi}(r y_{\text{fin}}) }{\sqrt{\N(r y_{\text{fin}})}} W_{\phi}(r y_{\infty}) \psi(rx)
,\qquad y\in \mathbb{A}^{\times},\ \ x \in \mathbb{A},\]
\[W_\phi:=W_{(R^{(\mt)})^{-1}\phi}\quad\text{and}\quad
\lambda^{(\mt)}_{\pi}(\mathfrak{m}) := \sum_{\mathfrak{s}
\mid\gcd(\mt,\mathfrak{m})} \alpha_{\mt, \mathfrak{s}}
(\N\mathfrak{s})^{1/2} \lambda_{\pi}(\mathfrak{m}\mathfrak{s}^{-1}).\]
The constants $\alpha_{\mt, \mathfrak{s}}\in\CC$ are explicit, although difficult to estimate in general.

Similarly, the Eisenstein spaces have an orthogonal decomposition
\[V_{\chi,\chi^{-1}}(\mc)=\bigoplus_{\mt\mid\mc\mc_\chi^{-2}}
R^{(\mt)} V_{\chi,\chi^{-1}}(\mc_\chi^2), \qquad \mc_\chi^2\mid\mc,\]
such that every $\phi \in R^{(\mt)}V_{\chi,\chi^{-1}}(\mc_\chi^2)$ has a Fourier expansion
\[\phi\left( \left(\begin{matrix} y & x\\ 0 &
1\end{matrix}\right)\right) = \phi_{\text{const}}(y) + \sum_{r \in
K^\times} \frac{\lambda^{(\mt)}_{\chi,\chi^{-1}}(r
y_{\text{fin}})}{\sqrt{\N (r y_{\text{fin}})}} W_{\phi}(r
y_{\infty}) \psi(rx),\qquad y\in \mathbb{A}^{\times},\ \ x \in
\mathbb{A},\]
\[\lambda^{(\mt)}_{\chi,\chi^{-1}}(\mathfrak{m})\ll_{K,\eps}
(\N\gcd(\mt,\mathfrak{m}))(\N\mathfrak{m})^\eps.\]

\medskip

\section{Lecture Four: Subconvexity of Twisted $L$-functions}

In this final lecture we highlight the main ideas in the proof of Theorem~\ref{thm1}, following closely the original source \cite{BH2}. The proof is inspired by and builds on important earlier work by several researchers: see \cite{BH3,BH2} for references and history.
For the sake of readability, we omit some technicalities, and we do not indicate the dependence of implied constants on $\pi$, $\chi_\infty$, $K$, $\eps$.

By an approximate functional equation it suffices to bound finite sums
\begin{equation}\label{eq13}
\mathcal{L}_{\chi_\text{fin}}:=
\sum_{0 << r \in \my}\frac{\lambda_{\pi}(r\my^{-1})\tilde\chi_\text{fin}(r)}{\sqrt{\N(r\my^{-1})}}W(r),
\end{equation}
where $W:\Kplus\to\CC$ is a smooth weight function of compact support cutting off at about $\N r\approx(\N\mq)^{1+\eps}$, and $\my$ represents a narrow ideal class. Recall that $\tilde\chi_\text{fin}:(\mo/\mq)^\times\to S^1$ is determined by the Hecke character $\chi:\aA^\times\to S^1$. We shall estimate \eqref{eq13} by a general principle of analytic number theory: sums in a harmonic family are easier to bound together than individually. To illustrate this point, let us assume that we are back in school, and we need to prove the inequality
\[|\sin x+\cos x|\leq\sqrt{2}.\]
Of course, we can accomplish this task in many ways, but a particularly nice way is to ``invent the harmonic complement'' $|\sin x-\cos x|$ and observe the identity
\[(\sin x+\cos x)^2+(\sin x-\cos x)^2=2.\]

In order to estimate \eqref{eq13}, we consider \emph{all the sums}
\[\mathcal{L}_{\xi}:= \sum_{0 << r \in \my} \frac{\lambda_{\pi}(r\my^{-1})\xi(r)}{\sqrt{\N(r\my^{-1})}}W(r),\]
where $\xi$ is any character of $(\mo/\mq)^{\times}$, and then the \emph{amplified second moment}
\[\cS:=\sum_{\xi \in \widehat{(\mo/\mq)^{\times}}}\left|\sum_{\ell} \xi(\ell)\ov{\tilde\chi_\text{fin}(\ell)}\right|^2 \bigl|\mathcal{L}_{\xi}\bigr|^2,\]
where $\ell\in\mo$ runs through certain elements of norm $\N\ell\in[L,2L]$ generating prime ideals $(\ell)\nmid\mq$. The quantity $L>1$ is the \emph{amplifier length}, and the previous sum is
\[\cS\gg(\N\mq)^{-\eps}L^2\bigl|\mathcal{L}_{\tilde\chi_\text{fin}}\bigr|^2\]
by positivity. Applying Plancherel and some easy estimates for the diagonal contribution, we arrive at
\begin{equation}\label{eq14}
\frac{|\mathcal{L}_{\tilde\chi_\text{fin}}|^2}{(\N\mq)^{1+\eps}}\ll \frac{1}{L}+
\sum_{0\neq q \in\mq\my\cap\mathcal{B}}\ \sum_{\substack{\ell_1r_1 - \ell_2r_2 =q \\ 0\neq r_1, r_2 \in \my}}
\frac{\lambda_{\pi}(r_1\my^{-1})\ov{\lambda_\pi(r_2\my^{-1})}}{\sqrt{\N (r_1r_2\my^{-2})}}W(r_1)\ov{W(r_2)},
\end{equation}
where $(\ell_1)$, $(\ell_2)$ are some prime ideals of norm about $L$, and $\mathcal{B}\subset K_\infty$ is some ball of volume at most $L(\N\mq)^{1+\eps}$ centered at the origin.

To handle the \emph{shifted convolution sum} inside \eqref{eq14}, we write it as
\begin{equation}\label{eq15}
\sum_{\substack{r_1-r_2=q \\ r_1, r_2 \in K^\times}}
\frac{\lambda_{\pi}(r_1\ell_1^{-1}\my^{-1})\ov{\lambda_\pi(r_2\ell_2^{-1}\my^{-1})}}
{\sqrt{\N (r_1\ell_1^{-1}r_2\ell_2^{-1}\my^{-2})}}W_1(r_1)\ov{W_2(r_2)}.
\end{equation}
The weight functions $W_1,W_2:K_\infty^\times\to\CC$ are nice, hence they determine vectors $\phi_1,\phi_2\in V_\pi(\mc_\pi)$ such that
\[\phi_j\left(\left(\begin{matrix} y & x\\ 0 &1\end{matrix}\right)\right)=
\sum_{r \in K^\times}\frac{\lambda_{\pi}(r y_{\text{fin}})}{\sqrt{\N (r y_{\text{fin}})}}
W_j(r y_{\infty}) \psi(rx),\qquad y\in \mathbb{A}^{\times},\ \
x \in \mathbb{A}.\]
Let us fix $y\in \mathbb{A}^{\times}$ such that $y_\infty=(1,\dots,1)$ and $(y_{\text{fin}})=\my^{-1}$. Let us also define $\Phi\in L^2(\GL_2(K)\backslash \GL_2(\mathbb{A}))$ via
\[\Phi(g):=\phi_1\left(g\begin{pmatrix}\ell_1^{-1}&0\\0& 1\end{pmatrix}\right)
\ov{\phi_2}\left(g\begin{pmatrix}\ell_2^{-1}&0\\0& 1\end{pmatrix}\right),\qquad g\in\GL_2(\aA).\]
It is straightforward to check that $\Phi$ is fixed by the right action of $\KK(\mc)$ for
\[\mc:=\mc_\pi\lcm((\ell_1),(\ell_2)).\]
Moreover,
\[\Phi\left(\left(\begin{matrix} y & x\\ 0 &1\end{matrix}\right)\right)=
\sum_{r_1,r_2\in K^\times}\frac{\lambda_{\pi}(r_1\ell_1^{-1}\my^{-1})\ov{\lambda_\pi(r_2\ell_2^{-1}\my^{-1})}}
{\sqrt{\N (r_1\ell_1^{-1}r_2\ell_2^{-1}\my^{-2})}}W_1(r_1)\ov{W_2(r_2)}\psi((r_1-r_2)x),\]
whence \eqref{eq15} really equals
\[\int_{K\backslash\mathbb{A}} \Phi\left(\begin{pmatrix} y & x\\ 0 &1\end{pmatrix}\right) \psi(-qx) \,dx.\]
We decompose $\Phi$ in the level $\mc$ spectrum of $L^2(\GL_2(K)\backslash \GL_2(\mathbb{A}))$, following the discussion of the third lecture. We obtain
\[\Phi = \Phi_\text{sp}+\int_{(\mc)}\sum_{\mt\mid\mc\mc_\varpi^{-1}}\Phi_{\varpi,\mt}\,d\varpi,\]
where $\Phi_\text{sp}$ lies in the span of the functions $g\mapsto\chi(\det g)$ for quadratic Hecke characters $\chi$, the representations $V_\varpi$ run through the cuspidal spaces $V_\pi$ and the Eisenstein spaces $V_{\chi,\chi^{-1}}$ of conductor $\mc_\varpi\mid\mc$, and
\[\Phi_{\varpi,\mt}\in R^{(\mt)}V_\varpi(\mc_\varpi),\qquad \mt\mid\mc\mc_\varpi^{-1}.\]
Upon defining $W_{\varpi,\mt}:=W_{\Phi_{\varpi,\mt}}$, we can now rewrite \eqref{eq14} as
\begin{equation}\label{eq16}
\frac{|\mathcal{L}_{\tilde\chi_\text{fin}}|^2}{(\N\mq)^{1+\eps}}\ll \frac{1}{L}+
\sum_{0\neq q \in\mq\my\cap\mathcal{B}}\ \int_{(\mc)} \sum_{\mt\mid\mc\mc_{\varpi}^{-1}}
\frac{\lambda_{\varpi}^{(\mt )}(q\my^{-1})}{\sqrt{\N(q\my^{-1})}} W_{\varpi,\mt}(q)\,d\varpi.
\end{equation}
The advantage of this expression is that it is not quadratic but linear in the Hecke eigenvalues. The price to pay is the spectral averaging in $\varpi$.

The Whittaker functions $W_{\varpi,\mt}:K_\infty^\times\to\CC$ depend on the original weight function $W:\Kplus\to\CC$ included in \eqref{eq13}. In order to proceed further, one needs to understand the ``size'' of these functions, at least on average over the spectrum. We carry this out via Sobolev type norms, the process schematized as follows:
\[\|W\|_\ast\nyil\|\phi_1\|_\ast\|\phi_2\|_\ast\nyil\|\Phi\|_\ast\nyil\|\Phi_{\varpi,\mt}\|_\ast\nyil\|W_{\varpi,\mt}\|_\ast\]
For example, over the Eisenstein spectrum we can readily derive the bound
\[\int\limits_{\varpi\in\mathcal{E}(\mc)} \sum_{\mt \mid \mc\mc_{\varpi}^{-1}}
\left|W_{\varpi, \mt }(y)\right| \,d\varpi \ll(\N(\ell_1\ell_2))^{\eps},\]
using that the number of cusps in level $\mc$ is $\ll(\N(\ell_1\ell_2))^{\eps}$, because $\mc\mc_\pi^{-1}$ is square-free.
This implies that the contribution of the Eisenstein spectrum in \eqref{eq16} is
\[\ll\frac{1}{\N\mq}\sqrt{L(\N\mq)^{1+\eps}}(\N(\ell_1\ell_2))^{\eps}\ll\frac{L^\frac{1}{2}}{(\N\mq)^{\frac{1}{2}-\eps}}.\]
Bounding the cuspidal contribution in \eqref{eq16} is much harder, but at least we can initially restrict to very small spectral parameters, namely $\N\tilde{\lambda}_{\varpi}\leq(\N\mq)^\eps$, thanks to the bound (valid for any $A>0$)
\[\int\limits_{\varpi\in\mathcal{C}(\mc)}(\N\tilde{\lambda}_{\varpi})^A
\sum_{\mt \mid \mc\mc_{\varpi}^{-1}}\left|W_{\varpi, \mt}(y)\right| \,d\varpi \ll_A|\N(\ell_1\ell_2)|^{\frac{1}{2}+\eps}.\]
Using this observation we find that the contribution of the cuspidal spectrum in \eqref{eq16} is, for some values $f(\ma)\ll (\N\mq)^\eps$,
\[\ll(\N\mq)^\eps\left(\sum_{\varpi\in\mathcal{C}(\mc,\eps)}\sum_{\mt\mid\mc\mc_{\varpi}^{-1}} \left|\sum_{\N\mathfrak{m} \ll L(N\mq)^\eps}
\frac{\lambda_{\varpi}^{(\mt )}(\mathfrak{m}\mq)}{\sqrt{\N(\mathfrak{m}\mq)}}f(\mathfrak{m}\mq)\right|^2\right)^{1/2}.\]
Here we essentially factor out $\lambda_{\varpi}(\mq )$ and bound it individually by $(\N\mq)^{\theta+\eps}$, which is why the parameter $\theta$ appears in Theorems~\ref{thm1} and~\ref{thm2}. Now we arrive at the endgame. We majorize the remaining sum by a smooth spectral sum involving an analogous Eisenstein contribution and rapidly decaying spectral weights. We open the square and apply a variant of the Bruggeman--Kuznetsov formula originally developed by Venkatesh~\cite{Ve2} and further extended by Maga~\cite{Ma}. Finally we apply familiar bounds for the resulting Kloosterman sums and Bessel transforms to infer that the contribution of the cuspidal spectrum in \eqref{eq16} is
\[\ll(\N\mq)^\eps\Biggl(\underbrace{\frac{L^2}{(\N\mq)^{1-2\theta-\eps}}}_{\text{diagonal}}
+\underbrace{\frac{L^{3/2}}{(\N\mq)^{1-2\theta-\eps}}}_{\text{off-diagonal}}\Biggr)^{1/2}\ll\frac{L}{(\N\mq)^{\frac{1}{2}-\theta-\eps}}.\]

\smallskip

Collecting terms, we deduce from \eqref{eq16} that
\[\frac{|\mathcal{L}_{\tilde\chi_\text{fin}}|^2}{(\N\mq)^{1+\eps}}\ll \frac{1}{L}+\frac{L}{(\N\mq)^{\frac{1}{2}-\theta}}.\]
The right hand side is smallest when $L:=(\N\mq)^{\frac{1-2\theta}{4}}$, in which case we get
\[\mathcal{L}_{\chi_\text{fin}}\ll(\N\mq)^{\frac{3+2\theta}{8}+\eps}.\]
This concludes the proof of Theorem~\ref{thm1}.

\medskip

\begin{acknowledgements}
I thank my friends Valentin Blomer and Jianya Liu for making these notes possible. Valentin and I worked on the results between 2007 and 2009, and we learnt a lot from each other. Jianya invited me to lecture at a great summer school in 2010, and he showed infinite patience in waiting for my written account. Finally, Rainer Schulze-Pillot explained to me some subtleties of quadratic forms over number fields, for which I am most grateful.
\end{acknowledgements}

\end{document}